\documentclass[11pt,amstex]{article}
\usepackage{amsmath}
\usepackage{amsfonts}
\usepackage{amscd}
\usepackage[dvips]{color}
\def\a{\alpha}
\def\b{\beta}
\def\ga{\gamma}

\def\la{\lambda}

\def\phi{\varphi}

\def\be{\begin{equation}}
\def\ee{\end{equation}}
\def\bear{\begin{eqnarray}}
\def\eear{\end{eqnarray}}
\def\best{\begin{eqnarray*}}
\def\eest{\end{eqnarray*}}
\def\pf{{\bf Proof}: }
\newtheorem{thm}{Theorem}

\newtheorem{theorem}{Theorem}[section]
\newtheorem{prop}[theorem]{Proposition}

\newtheorem{lemma}[theorem]{Lemma}
\newtheorem{cor}[theorem]{Corollary}

\newtheorem{remark}[theorem]{Remark}
\newenvironment{rem}{\begin{remark}\rm}{\end{remark}}
\newtheorem{example}[theorem]{Example}
\newenvironment{ex}{\begin{example}\rm}{\end{example}}

\def\non{\noindent}
\def\pf{\non {\bf Proof. }}
\def\qed{\nopagebreak \hskip .1in { $\Box$ }\penalty10000 %
\hskip\parfillskip \par  }

\def\ti{\times}
\def\del{\overline \partial}
\def\bd{\partial}

\def\ker{\mbox{ker\,}}

\def\ov#1{\overline{#1}}

\def\R{{\mathbb R}}
\def\P{{\mathbb P}}
\def\Q{{\mathbb Q}}
\def\cx{{\mathbb C}}

\def\F{{\cal F}}

\def\M{{\cal M}}

\newcommand{\CM}{\overline{{\cal M}}}

\title{\bf Holomorphic 2-forms and Vanishing Theorems for
           Gromov-Witten Invariants} \vskip.2in
\author{ Junho Lee }
\date{\empty}
\addtocounter{section}{0}
\begin{document}
\maketitle
\vskip.1in
\begin{abstract}
On a  compact K\"{a}hler manifold $X$ with  a holomorphic 2-form
$\a$, there is an almost complex structure  associated with $\a$. We
show how this implies vanishing theorems for the Gromov-Witten
invariants of $X$. This  extends  the approach, used in \cite{lp}
for K\"{a}hler surfaces, to higher dimensions.
\end{abstract}
\vskip.3in


Let $X$ be a K\"{a}hler surface with a non-zero holomorphic 2-form
$\a$. Then $\a$ is a section of  the canonical bundle and its   zero
locus $Z_\a$, with multiplicity, is  a canonical divisor.  We showed
in  \cite{l1} that the real 2-form $\mbox{Re}(\a)$ determines a
(non-integrable) almost complex structure $J_\a$ that has the
following remarkable  ``Image Localization Property''\,: {\em if a
$J_\a$-holomorphic map $f:C\to X$ represents a non-zero (1,1) class,
then $f$ is  in fact holomorphic and its image lies in $Z_\a$.}  As
shown in \cite{lp}, this property  together with Gromov Convergence
Theorem leads to\,:
\begin{thm}[\cite{lp}]\label{T:lp}
Let $X$ be a K\"{a}hler surface with a non-zero holomorphic 2-form
$\a$. Then, any class $A$ with non-trivial Gromov-Witten invariant
$GW_{g,k}(X,A)$  is represented by a stable holomorphic map $f:C\to
X$ whose image lies in the canonical divisor $Z_\a$.
\end{thm}

This paper  extends Theorem~\ref{T:lp}  to higher dimensions. The
principle is the same: perturbing the K\"{a}hler structure to a
non-integrable almost complex structure $J_\a$ forces the
holomorphic maps to satisfy  certain geometric conditions determined by $\a$.
This gives constraints on the Gromov-Witten invariants.

Specifically, let $X$ be a compact K\"{a}hler manifold with a non-zero holomorphic
2-form $\a$. Then the real part of $\a$ defines an endomorphism
$K_\a$ of $TX$ and  an almost complex structure $J_\a$, just as in
the surface case (see (\ref{al-cx-str}) and (\ref{end})).
These geometric structures lead, naturally and easily, to our main theorem\,:

\begin{thm}\label{ext}
Let $X$ be a  compact K\"{a}hler manifold  with a non-zero
holomorphic 2-form $\a$. Then any class $A$ with non-trivial
Gromov-Witten invariant $GW_{g,k}(X,A)$ is represented by a stable
holomorphic map $f:C\to X$ satisfying the equation $ K_\a df=0.$
\end{thm}

This theorem  follows from Theorem~\ref{T:Main}, which is more suitable for applications.  It generalizes
Theorem~\ref{T:lp} since when $X$ is a surface  the kernel
of the endomorphism $K_\a$ is trivial  on $X\setminus Z_\a$ (see Example~\ref{E:lp}).
The equation $ K_\a df=0$ is a geometric fact about holomorphic maps that directly implies numerous
vanishing results about Gromov-Witten invariants (see Section 3).

 Section 1 briefly describes  Gromov-Witten
invariants and states a vanishing principle for them. Section 2 contains
the definition of the  almost complex structures $J_\a$ and some of the
consequences of that definition.  In Section 3
we apply  a stronger version of Theorem 2, which directly follows from
properties of $J_\a$, to show various vanishing results for
 Gromov-Witten invariants.

\medskip
\noindent {\bf Acknowledgments.} I am very grateful to T. Parker for
useful discussions and encouraging conversations.
I am also very thankful to the referee for several helpful comments and corrections.

\vskip1cm

\setcounter{equation}{0}
\section{Gromov-Witten Invariants}
\label{section1}
\bigskip

The aim of this section is to give a brief description of the  Gromov-Witten
invariants and to set up notations for them.
Let $(X,\omega)$ be a compact symplectic $2n$-dimensional manifold with
an  $\omega$-tamed almost complex structure $J$, i.e.,
$\omega(u,Ju)>0$. A $J$-holomorphic map
$f:(C,j)\to X$ from a  (connected) marked nodal curve is
{\em stable}  if  its automorphism group is finite\,(cf. \cite{hz}).
Denote by
\begin{equation*}\label{sms}
 \CM_{g,k}(X,A,J)
\end{equation*}
the moduli space of stable $J$-holomorphic maps from marked
nodal curves of (arithmetic) genus $g$ with $k$ marked points that
represent the homology class $A\in H_2(X)$.
This moduli space carries a (virtual) fundamental homology class
of real dimension
\begin{equation}\label{dim}
 2\,\left[\,c_1(TX)\cdot A+(n-3)(1-g)+k\,\right]
\end{equation}
(cf. \cite{lt})
whose push-forward under the map
\begin{equation*}\label{st-ev}
 st \times ev:\CM_{g,k}(X,A,J) \to
\CM_{g,k} \times X^k
\end{equation*}
defined by stabilization and evaluation at the marked points is the
Gromov-Witten invariant
\begin{equation}\label{GW}
GW_{g,k}(X,A) \ \in\  H_*(\ov{\M}_{g,k}\ti X^k;\Q).
\end{equation}
This is equivalent to the collection of  ``GW numbers''
\begin{equation*}\label{GWn}
GW_{g,k}(X,A) (\mu; \gamma_1,\dots, \gamma_k)
\end{equation*}
obtained by evaluating the homology class (\ref{GW}) on  the
cohomology classes Poincar\'{e} dual to $\mu\in H_*(\CM_{g,k})$ and
$\gamma_j\in H_*(X)$ whose total degree is the dimension
(\ref{dim}). Standard cobordism arguments  then show that these are
independent of the choice of  $J$, and depend only on
the deformation class of the symplectic form $\omega$.

Our subsequent discussions are based on the following
vanishing principle for GW invariants.

\begin{prop}\label{vanishing}
Fix a compact symplectic manifold $(X,\omega)$. Suppose
$$ GW_{g,k}(X,A)(\mu;\ga_1,\cdots,\ga_k)\  \ne\  0.$$
Then, for any $\omega$-tamed almost complex structure $J$ and for
any geometric representatives $M\subset \CM_{g,k}$ and
$\Gamma_j\subset X $ of  classes $\mu\in H_*(\CM_{g,k})$ and
$\gamma_j\in H_*(X)$ there exists a stable $J$-holomorphic map
$f:(C,x_1,\cdots,x_k)\to X$ representing class $A$ with $st(C)\in M$
and $f(x_j)\in \Gamma_j$.
\end{prop}

The proof is straightforward (cf. \cite{lp}). For convenience, we
will assemble all GW invariants for a class $A$ into a single
invariant by introducing a variable $\la$ to keep track of the
genus. The GW series of $(X,\omega)$ for a class $A$ is then the
formal power series
\begin{equation*}\label{formal-s}
 GW_A(X)\ =\ \sum_{g,k}\,\frac{1}{k!}\,GW_{g,k}(X,A)\,\la^{g}.
\end{equation*}

\vskip1cm
\setcounter{equation}{0}
\section{The Almost Complex Structures $J_\a$ }
\label{section2}
\bigskip

Let $(X,\omega)$ be a compact symplectic manifold with an
$\omega$-compatible almost complex structure $J$, namely  $\langle
u,v\rangle=\omega(u,Jv)$ is a Riemannian metric. A 2-form $\a$ is
then called   {\em $J$-anti-invariant} if $\a(Ju,Jv)  =  -\a(u,v)$.
As observed in \cite{l1}, each $J$-anti-invariant 2-form $\a$
induces an almost complex structure
\begin{equation}\label{al-cx-str}
 J_\a\  = \ (Id+JK_\a)^{-1}J(Id+JK_\a)
\end{equation}
where $K_\a$ is an endomorphism of $TX$ defined by the equation
\begin{equation}\label{end}
 \langle u,K_\a v\rangle\    =\   \a(u,v).
\end{equation}
Such endomorphisms $K_\a$ are skew-adjoint and anti-commute with
$J$. It follows that $Id+JK_\a$ is invertible and hence
(\ref{al-cx-str}) is well-defined. A simple computation then shows
that for any $C^1$ map $f:(C,j)\to X$,
\begin{equation}\label{equiv1}
\del_{J_a}f = 0\ \ \ \Longleftrightarrow\ \ \
\del_J f   =    K_\a\bd_J f j
\end{equation}
where
$$
\del_J f \ =\ \tfrac12(df+Jdf)\,,\ \ \ \ \ \
\bd_J f \ =\  \tfrac12(df-Jdf j).
$$
(\ref{equiv1}) implies that every
$J$-holomorphic map $f$ satisfying $K_\a df=0$ is also
$J_\a$-holomorphic.
One can also show that if $f$ is $J_\a$-holomorphic then
\begin{equation}\label{int-equiv}
\int_{C}|\del_J f|^2\ =\ \int_{C}|K_\a\bd_J f|^2\ =\
\int_{C}f^*(\a)
\end{equation}
(cf. \cite{l1}). This integral vanishes when $\a$ is closed and $\a(A)=0$ where $A$
is the class represented by $f$. In this case, the given $J_\a$-holomorphic map
$f$ is $J$-holomorphic ($\del_J f=0$) and satisfies  $K_\a
df=K_\a\bd_J f=0$. Therefore, when $\a$ is closed and $\a(A)=0$,
a map $f$ representing the class $A$ is $J_\a$-holomorphic if and
only if $f$ is $J$-holomorphic and satisfies the equation $K_\a
df=0$. Combined with Proposition~\ref{vanishing}, these observations
lead to\,:

\begin{prop}\label{P:Main}
Let $(X,\omega)$ be a compact symplectic manifold with an
$\omega$-compatible $J$ and with a closed $J$-anti-invariant 2-form
$\a$. Then, for any class $A$ with $GW_{g,k}(X,A)\ne 0$ we have
\begin{equation*}\label{T:msp}
 \CM_{g,k}(X,A,J_\a)\ =\ \{\ f\in\CM_{g,k}(X,A,J)\ |\ K_\a df = 0\ \}.
\end{equation*}
In particular, this space is not empty.
\end{prop}

\pf By the above discussion, it suffices to show that $\a(A)=0$ and
$\CM_{g,k}(X,A,J_\a) \ne \emptyset$. Proposition~\ref{vanishing}
shows that there is a $J$-holomorphic map $h:(D,j)\to X$
representing the class $A$. Fix a point $p\in D$ and choose an
orthonormal basis $\{e_1,e_2=je_1\}$ of $T_pD$. Then,
\begin{equation*}
h^*\a(e_1,e_2)\ =\ \a(h_*e_1,h_*je_1)\ =\ \a(h_*e_1,Jh_*e_1).
\end{equation*}
Since $\a$ is $J$-anti-invariant, this vanishes and hence
$\a(A)=\int_D h^*(\a)=0$. On the other hand, for any sufficiently
small $t>0$ the almost complex structure $J_{t\a}$ is
$\omega$-tamed since $\omega$-tamed is an open condition.
Proposition~\ref{vanishing} then asserts that there exists a
$J_{t\a}$-holomorphic map $f$ representing the class $A$. By
(\ref{int-equiv}) and the fact $K_{t\a}=tK_\a$, this map $f$ is
$J$-holomorphic  and satisfies $K_\a df= 0$. Thus, $f$ is also
$J_\a$-holomorphic by (\ref{equiv1}). \qed

\bigskip

Below, we will show some basic properties of the zero locus $Z_\a$
of $\a$ and $\mbox{ker\,}K_\a$, which will be frequently used in
subsequent arguments. One can use $J$ to decompose
$\Omega^2(X)\otimes \cx$ as
\begin{equation*}
 \Omega^2(X)\otimes \cx\ =\
 \Omega^{2,0}_J(X)\oplus\Omega^{1,1}_J(X)\oplus\Omega^{0,2}_J(X).
\end{equation*}
Every $J$-anti-invariant 2-form $\a$ then can  be written as
$\a=\b+\overline{\b}$ for some $\b\in \Omega^{2,0}_J(X)$.
The following lemma simply follows from the definitions and the
properties of $K_\a$.

\begin{lemma}\label{A1}
Let $\mbox{dim\,}X=2n$, and $\a$ and $\b$ be as above. Then,
\begin{enumerate}
\item[(a)] $\a$ and $\b$ have the same zero locus, 
\item[(b)] if $n$ is odd then $\a^n=0$, and
           if $n=2m$ then $\a^n=c\,\b^m\wedge\overline{\b}^m$
           where $c= {n\choose m}$,
\item[(c)] the (real) dimension of $\mbox{ker}\,K_{\a}$ is at most $2n-4$
           at every point in   $X\setminus Z_\a$,
\item[(d)] $u\in \mbox{ker\,}K_\a$ if and only if
           $\a(u,w)=0$ for any $w$. Thus,
           $\mbox{ker}\,K_{\a}$ is trivial if and only if
           $\a$ is non-degenerate.
\end{enumerate}
\end{lemma}

A foliation $\F$ of dimension $m$ on $n$-dimensional manifold $M$ is
a decomposition $\F=(L_i)_{i\in I}$ of $M$ into pairwise disjoint
connected subsets $L_i$, which are called leaves of the foliation
$\F$, with the following property\,: for each $p\in M$ there exists
a foliation chart $\varphi:U\to W_1\times W_2$, where $W_1$ and
$W_2$ are open disks in $\R^m$ and $\R^{n-m}$ respectively, such
that for each point $q\in W_2$ the preimage
$\varphi^{-1}(W_1\times\{q\})$ is a connected component of $U\cap
L_i$ for some leaf $L_i$. We refer to \cite{cn} and \cite{hol} for more details on
foliations.

\begin{lemma}\label{six}
Let $(X,\omega)$ be a six dimensional symplectic manifold with $\omega$-compatible
$J$. If $\a$ is a closed $J$-anti-invariant 2-form, then
$\mbox{ker\,}K_\a$ gives a foliation on $X\setminus Z_\a$ of (real) dimension
two whose leaves are all $J$-invariant.
\end{lemma}

\pf Since  $K_\a$ is anti-commute with $J$,
Lemma~\ref{A1}\,c implies that on  $X\setminus Z_\a$
the dimension of $\ker K_\a$ is  two. On the other hand,
$d\a(u,v,w)=0$  gives
\begin{equation*}
L_u(\a(v,w)) - L_v(\a(u,w)) + L_w(\a(u,v))
           - \a([u,v],w) + \a([u,w],v) - \a([v,w],u)=0
\end{equation*}
where $L$ denotes the Lie derivative. This together with
Lemma~\ref{A1}\,d imply that if $u,v\in \mbox{ker\,}K_\a$ then
$[u,v]\in \mbox{ker\,}K_\a$. Therefore, by Frobenius Theorem
$\mbox{ker\,}K_\a$ gives a foliation on $X\setminus Z_\a$ of dimension two.
Since  $K_\a$ is anti-commute with $J$, every leaf is $J$-invariant. \qed

\vskip1cm
\setcounter{equation}{0}
\section{Vanishing Results}
\label{section3}
\bigskip

Let $(X,J)$ be a compact K\"{a}hler manifold with a non-zero
holomorphic 2-form $\a$. By the Hodge Theorem $\a$ is closed and hence
its real part $\mbox{Re}(\a)$ is also closed. Moreover, the real
2-form $\mbox{Re}(\a)$ is $J$-anti-invariant and its zero locus  is
$Z_\a$ by Lemma~\ref{A1}\,a.
Throughout this section, we will denote by $K_\a$ the endomorphism of $TX$
defined by $\mbox{Re}(\a)$ as in (\ref{end}).

\medskip
A holomorphic 2-form $\a$ is called {\em non-degenerate} if
$\mbox{Re}(\a)$ is non-degenerate, or equivalently $\ker K_\a$ is trivial.
The following theorem directly follows from
Proposition~\ref{vanishing} and Proposition~\ref{P:Main}.

\begin{theorem}\label{T:Main}
Fix a compact K\"{a}hler manifold $X$ with a non-zero holomorphic
2-form $\a$. If for a non-zero class $A$
\begin{equation*}
 GW_{g,k}(X,A)(\mu;\ga_1,\cdots,\ga_k)\ \ne\  0
\end{equation*}
then for any  geometric representatives $M\subset \CM_{g,k}$ and
$\Gamma_j\subset X $ of classes $\mu\in H_*(\CM_{g,k})$ and
$\gamma_j\in H_*(X)$ there exists a stable holomorphic map
$f:(C,x_1,\cdots,x_k)\to X$ representing  the class $A$ with
$st(C)\in M$ and $f(x_j)\in \Gamma_j$ and satisfying the equation
$
K_\a\,df  =  0.
$
Consequently,  if $\a$ is
non-degenerate on an open set $U\subset X$ then the image of $f$ lies in
$X\setminus U$.
\end{theorem}

Using this Theorem, one can obtain various vanishing results
about GW invariants.

\begin{ex}
Given a compact hyperk\"{a}hler manifold $X$ of (complex) dimension
$2m$, there exists a holomorphic symplectic 2-form $\a$, i.e. $\a^m$
is nowhere vanishing (cf. \cite{bdl}). The 2-form $\a$ is  non-degenerate
on $X$ and hence Theorem~\ref{T:Main} implies that the series   $GW_A(X)$
vanishes unless  $A = 0$.
\end{ex}

\begin{ex}\label{prod}
Let $X=E_1\times \cdots \times E_n$  where each $E_i$ is an elliptic
curve and $n\geq 2$. For $i\ne j$, denote by $\a_{ij}$ the pull-back
2-form $\pi_i^*(\la_i)\wedge\pi_j^*(\la_j)$ where $\pi_i:X\to E_i$
is the $i$-th projection and $\la_i$ is a nowhere vanishing
holomorphic 1-form on $E_i$.
Now, suppose $GW_A(X)\ne 0$. Theorem~\ref{T:Main} then shows that
there is a holomorphic map $f:C\to X$ representing the class $A$ with $K_{\a_{ij}} df =0$.
Since $\a_{ij}$ has no zeros and $\ker K_{\a_{ij}}$ consists of vectors tangent to
fibers of the projection $\pi_i\times \pi_j:X\to E_i\times E_j$, we have
$(\pi_i\times \pi_j)_*df=0$ for each $i\ne j$. This implies  $A=0$.
The same arguments also apply to show that
 when $X=X_1\times\cdots \times X_n$ where each $X_i$ is a hyperk\"{a}hler manifold
or a complex torus of (complex) dimension at least two the series
$GW_A(X)$ vanishes unless $A= 0$.
\end{ex}

\begin{rem} There are well-known proofs for the above two examples (cf. \cite{bl}).
For instance, if $X$ is a compact hyperk\"{a}hler manifold
then every K\"{a}hler structure $J$
in the twistor family is deformation equivalent to $-J$ through K\"{a}hler structures
(cf. \cite{bdl}).
This directly implies $GW_A(X)=0$ unless $A=0$.
The product formula of \cite{b} for GW invariants  applies to
give the same vanishing results as in Example~\ref{prod}.
\end{rem}

The following example appears in \cite{lp}.

\begin{ex}(\cite{lp})\label{E:lp}
Let $X$ be a K\"{a}hler surface with a non-zero holomorphic 2-form
$\a$. Then, $\a$ is non-degenerate on $X\setminus Z_\a$ by
Lemma~\ref{A1}\,c,d. Note that since $\a$ is a section of the
canonical bundle the zero locus $Z_\a$ is a support of a canonical
divisor. Theorem~\ref{T:Main} thus shows that for any non-zero class
$A$ and for any genus $g$
\begin{equation}\label{pt-c}
  GW_{g,k}(X,A)(\,\cdot\,;\ga,\cdots)\ =\ 0
\end{equation}
where $\ga$ lies in $H_{i}(X)$ for $i=0,1$.
On the other hand, if $X$ is minimal surface of general type then
every canonical divisor is  connected (cf. \cite{bhpv}). We further
assume that the zero locus $Z_\a$ is a smooth (reduced) canonical
divisor. Then, any holomorphic map $f$ whose image lies in $Z_\a$
represents a (non-negative) multiple of the canonical class $K$.
Therefore, Theorem~\ref{T:Main} implies that the series $GW_A(X)$
vanishes unless $A=mK$ for some non-negative integer $m$.
\end{ex}

The following example extends both the vanishing result (\ref{pt-c}) and
  Theorem~\ref{T:lp} of
the introduction to  K\"{a}hler manifolds of even complex dimension.   It is an immediate consequence of Theorem~\ref{T:Main}.

\begin{ex}
Fix a compact K\"{a}hler manifold $X$ of complex dimension $2m$
with a holomorphic 2-form $\a$. If $\a^m$ is not identically zero,
then the zero locus $Z_m$ of $\a^m$, with multiplicities, is a canonical divisor of $X$ and
$\a$ is non-degenerate on $X\setminus Z_m$.  Theorem~\ref{T:Main}  implies that:
\begin{enumerate}
\item[(a)] if $GW_{g,k}(X,A)\ne 0$ for a non-zero class $A$, then
           $A$ is represented by a stable holomorphic map $f:C\to X$ whose image
           lies in the canonical divisor $Z_m$, and
\item[(b)] for any non-zero class $A$ and for any genus $g$ we have
           $ GW_{g,k}(X,A)(\,\cdot\,;\ga,\cdots)    =   0$
           where $\ga$ lies in $H_{i}(X)$ for $i=0,1$.
\end{enumerate}
\end{ex}

\bigskip
\non
{\bf Compact K\"{a}hler Threefolds}
\medskip

\non
Let $X$ be a compact K\"{a}hler threefold with a non-zero holomorphic 2-form $\a$.
It then follows from Lemma~\ref{six} that $\ker K_\a$ induces a foliation on $X\setminus Z_\a$
of (real) dimension two. We will denote this foliation by $\F_\a$.

\begin{lemma}\label{fol}
Fix a compact K\"{a}hler threefold $X$ with a non-zero holomorphic 2-form $\a$.
If a (non-constant) stable holomorphic map $f:C\to X$ satisfies the equation
$K_\a df=0$, then
\begin{itemize}
\item the image of each  irreducible component of $C$ either lies in $Z_\a$
or  lies in one  leaf of  the foliation $\F_\a$ on $X\setminus Z_\a$ union  finitely many
points of $Z_\a$.
\end{itemize}
Consequently, if $\a$ has no zeros then the image of $f$ lies in one leaf of the foliation
$\F_\a$ on $X$.
\end{lemma}

\pf
Collapse all irreducible components of $C$ whose image is a point. The resulting map
still has the same image $f(C)$, so we can assume that the image of each irreducible
component is not a point. Fix an irreducible component $C_i$ of $C$ and
suppose $f(C_i) $ is not contained in $Z_\a$.
Then, the intersection $f(C_i)\cap Z_\a$ is finite since $f$ is
holomorphic and $Z_\a$ is an analytic subvariety. Denote by $D_i$
the set of critical points of $f$ in $C_i$. This set  $D_i$ is
finite and hence $C_i\setminus (D_i \cup f^{-1}(Z_\a))$ is open and
connected. Therefore, the equation $K_\a df =0$ asserts that
 $f(C_i\setminus D_i) \setminus Z_\a\subset L_i$ for some leaf $L_i$
of the foliation $\F_\a$ on $X\setminus Z_\a$. It then remains to
show that for each $p\in D_i$ either $f(p)\in Z_\a$ or $f(p)\in
L_i$. Suppose $f(p)$ does not lie in $Z_\a$. Let  $(U,\varphi)$ be
a foliation chart around $f(p)$, namely $U\subset X\setminus Z_\a$
is a neighborhood of $f(p)$ and $\varphi(U)=W_1\times W_2\subset
\R^{2}\times \R^{4}$, where  $W_1$ and $W_2$ are open disks in
$\R^{2}$ and $\R^{4}$ respectively, such that for each point
$t\in W_2$ the pre-image $\varphi^{-1}(W_1\times \{t\})$ is a
connected component of $U\cap L_t$ for some leaf $L_t$ of $\F_\a$.
Then, for any small neighborhood $V\subset C_i$ of $p$ there exists
a point $t_i\in W_2$ such that $\varphi\circ f(V\setminus
\{p\})\subset W_1\times\{t_i\}$. Consequently, we have $\varphi\circ
f(p)\in W_1\times\{t_i\} $. Since the pre-image
$\varphi^{-1}(W_1\times\{t_i\})$  is a connected component of $U\cap
L_i$, we have $f(p)\in L_i$ . \qed

\begin{ex}\label{pro-bdl}
Fix a surface of general type $S$ with a holomorphic 2-form $\ga$
whose zero locus is a smooth canonical divisor $D$. Let $\pi:X=\P(TS)\to S$ be the
projective bundle with a pull-back 2-form $\a=\pi^*\ga$. The zero
locus $Z_\a$ is then the ruled surface
$\pi^{-1}(D)\to D$ and every leave of the
foliation $\F_\a$ on $X\setminus \pi^{-1}(D)$ is a  fiber of
$\pi:X\to S$. Thus, Theorem~\ref{T:Main} and Lemma~\ref{fol} together imply that
$GW_A(X)= 0$ unless $A=aD_0+bF$ for some integers $a$ and $b$ where
$D_0$ is the section class of the ruled surface $\pi^{-1}(D)$ and $F$ is the fiber class of $X$.
\end{ex}

Now, suppose $X$ is a compact threefold with a holomorphic 2-form
$\a$ without zeros. The foliation $\F_\a$ is then a foliation on the whole $X$.
In fact, $\F_\a$ is a holomorphic foliation\,; in a holomorphic local coordinates, the
foliation $\F_\a$ is given locally by the holomorphic vector field
\begin{equation*}
 Y\ =\ f_{23}\tfrac{\bd}{\bd z_1}\   -\
 f_{13}\tfrac{\bd}{\bd z_2}\   +\
 f_{12}\tfrac{\bd}{\bd z_3}
\end{equation*}
where $\a=f_{12}\, dz_1\wedge dz_2 + f_{13}\, dz_1\wedge dz_3
+f_{23} \,dz_2\wedge dz_3$.
After suitable change of local coordinates, we can  write
$Y=\tfrac{\bd}{\bd z_1}$. Such local coordinates give the required holomorphic
foliation chart.

\begin{cor}\label{cor}
Let $X$ be    a compact K\"{a}hler threefold with a holomorphic 2-form
$\a$ without zeros. Suppose $X$ is not a $\P^1$-bundle over a K3 or an abelian surface.
Then  for any non-zero class $A$ the invariant
\begin{equation}\label{pt-3}
 GW_{g,k}(X,A)(\,\cdot\,;\ga,\cdots)
\end{equation}
vanishes if the genus $g$ is $0$ or if one constraint $\ga$ lies in $H_{i}(X)$ for $0\leq i \leq 3$.
\end{cor}

\pf Assume  that for some  $A\neq 0$
the invariant (\ref{pt-3}) is not zero with
either  $g=0$ or
$\ga\in H_i(X)$ for $0\leq i\leq 3$.
We will show that
$X$ is a $\P^1$-bundle over a K3 or an abelian surface.
Since every leaf of $\F_\a$ is a smooth connected holomorphic curve,
by Theorem~\ref{T:Main} and Lemma~\ref{fol}
there exists a stable holomorphic map
$f:C\to X$ representing the non-zero class $A$ with $f(C)=L$
for some leaf $L$ of $\F_\a$. The leaf $L$ is thus compact and
$A=m[L]$ for some integer $m\geq 1$.
If the genus $g$ is $0$  then obviously $L=\P^1$.

On the other hand, if $\ga$ lies in $H_i(X)$ for $0\leq i\leq 3$ and the invariant (\ref{pt-3}) is non-zero, then the formal dimension (\ref{dim})  of the moduli space $\M_{g,0}(X,A)$ is strictly positive, so $c_1(X)\cdot A \geq 2$.  But by
Theorem 2 of \cite{le},  the normal bundle $N$ to $L$ satisfies $c_1(N)=0$, so
\begin{equation*}
2\leq  c_1(X)A\ =\ mc_1(X)[L]\ =\ m(\,c_1(L)+c_1(N)\,)[L]\ =\ mc_1(L) .
\end{equation*}
 Hence
$c_1(L)=2$ and therefore $L=\P^1$ in this case also.

 Now, by the proof of Corollary 2.8 of \cite{ho},
the fact that one  leaf $L$ of $\F_\a$ is a rational curve $\P^1$ implies that
every leaf of $\F_\a$ is  rational.  It follows that
the leaf space $S=X/\F_\a$ is a (smooth) compact K\"{a}hler surface and
the quotient map $\pi:X\to S$ is holomorphic.
Consequently, $X$ is a $\P^1$-bundle over $S$ and $\a$ descends to a
holomorphic 2-form $\ga$ on $S$ with $\pi^*\ga=\a$. Since the
holomorphic 2-form $\ga$ has no zeros, $c_1(S)=0$ and hence $S$ is a
K3 or an abelian surface (cf. \cite{bhpv}).  \qed

\begin{rem}
Let $X$ be a projective threefold with a holomorphic 2-form $\a$.
Suppose  that for some non-zero class $A$ the invariant
$GW_{g,k}(X,A)(\,\cdot\,;\ga,\cdots)$ with the constraint $\ga\in H_i(X)$ for $0\leq i\leq 3$
does not vanish. Then the canonical divisor $K_X$ is not nef by a dimension count.
In this case, there is a contraction $\pi:X\to S$ of an extremal ray such that
$X$ is a $\P^1$-bundle over a surface $S$ and $\a$ descends to a holomorphic 2-form $\ga$
on $S$ with $\pi^*\ga=\a$ (see Section 3 of \cite{cp}). Consequently, if $\a$ has
isolated  zeros then $S$ is a K3 or an abelian surface and, in fact, $\a$ has no zeros. This observation motivated
Corollary~\ref{cor}.
\end{rem}

\vskip 1cm

\bigskip
\noindent {\em Mathematics Department, Michigan State University, East Lansing, MI 48824}

\smallskip

\noindent {\em E-mail addresses:}\ \ {\ttfamily leejunho\@@msu.edu}

\end{document}